\documentclass[12pt]{article}%
\usepackage{amsfonts}
\usepackage{sw20bams}
\usepackage{amsmath}
\usepackage{amssymb}
\usepackage{graphicx}%
\providecommand{\U}[1]{\protect\rule{.1in}{.1in}}
\begin{document}

\title{M/D/1 Queues with LIFO and SIRO Policies}
\author{Steven Finch}
\date{October 18, 2022}
\maketitle

\begin{abstract}
While symbolics for the equilibrium M/D/1-LIFO waiting time density are
completely known, corresponding numerics for M/D/1-SIRO are derived from
recursions due to Burke (1959). \ Implementing an inverse Laplace
transform-based approach for the latter remains unworkable.

\end{abstract}

\footnotetext{Copyright \copyright \ 2022 by Steven R. Finch. All rights
reserved.}We examined in \cite{Fi-que2} a first-in-first-out M/D/1 queue
alongside unlimited waiting space, where the input process is Poisson with
rate $\lambda$ and the service times are constant with value $1/\mu=a$. \ The
FIFO policy of sorting clients, which is plainly fair and has minimal waiting
time variance \cite{Ki1-que2}, is also known as FCFS: first-come-first-serve.

Keeping all other features the same, we wonder about the effect of replacing
FIFO by LIFO:\ last-in-first-out. This policy, which seems patently unfair and
has maximal waiting time variance \cite{Tam-que2}, is also known as LCFS:\ last-come-first-serve.

Intermediate to FIFO\ and LIFO\ is SIRO: serve-in-random-order. \ This policy
is variously known as ROS:\ random-order-of-service and
RSS:\ random-selection-for-service. \ Analysis of SIRO\ is more difficult than
that of LIFO, as will be seen.

\section{LIFO}

Let $W_{\text{que}}$ denote the waiting time in the queue (prior to service).
\ Under equilibrium (steady-state) conditions and traffic intensity (load)
$\rho=\lambda/\mu<1$, the probability density function $f(x)$ of
$W_{\text{que}}$ has Laplace transform \cite{Tij-que2, Coo-que2, Cn1-que2,
Wis-que2}
\[
F(s)=\lim_{\varepsilon\rightarrow0^{+}}%
{\displaystyle\int\limits_{-\varepsilon}^{\infty}}
\exp(-s\,x)f(x)dx=1-\rho+\frac{\lambda-\lambda\,\Theta(s)}{s+\lambda
-\lambda\,\Theta(s)}=F_{\text{alt}}(s)+1-\rho
\]
where%
\[%
\begin{array}
[c]{ccc}%
\theta(x)=%
{\displaystyle\sum\limits_{k=1}^{\infty}}
e^{-k\,\rho}\dfrac{(k\,\rho)^{k-1}}{k!}\delta(x-a\,k), &  & \Theta
(s)=-\dfrac{1}{\rho}\,\omega\left(  -\rho\exp\left(  -\rho-a\,s\right)
\right)  ,
\end{array}
\]
$\delta(x)$ is the Dirac delta and $\omega(s)$ is the principal branch of the
Lambert omega:
\[%
\begin{array}
[c]{ccccc}%
\omega(s)e^{\omega(s)}=s, &  & -1\leq\omega(x)\in\mathbb{R}\text{ \ }%
\forall\text{ }x\geq-1/e, &  & \exists\text{ branch cut for }x<-1/e.
\end{array}
\]
From%
\[
(1-\rho)s+\lambda\,\left[  1-\Theta(s)\right]  (1-\rho+1)=s\,F(s)+\lambda
\,F(s)\left[  1-\Theta(s)\right]
\]
we have%
\[
(1-\rho)s+\lambda\left(  2-\rho-F(s)\right)  \left[  1-\Theta(s)\right]
=s\,F(s),
\]
i.e.,%
\[
F(s)=1-\rho+\lambda\left(  2-\rho-F(s)\right)  \left[  \frac{1}{s}%
-\frac{\Theta(s)}{s}\right]
\]
hence%
\begin{align*}
f(x)  &  =(1-\rho)\delta(x)+\kappa+\lambda%
{\displaystyle\int\limits_{0}^{x}}
\left(  (2-\rho)\delta(t)-f(t)\right)  \left[  1-%
{\displaystyle\int\limits_{0}^{x-t}}
\theta(u)du\right]  dt\\
&  =(1-\rho)\delta(x)+\lambda+\lambda(2-\rho)\left[  1-%
{\displaystyle\int\limits_{0}^{x}}
\theta(u)du\right]  -\lambda%
{\displaystyle\int\limits_{0}^{x}}
f(t)\left[  1-%
{\displaystyle\int\limits_{0}^{x-t}}
\theta(u)du\right]  dt.
\end{align*}
The indicated condition $\kappa=\lambda$ is true by the initial value theorem
\cite{Sch-que2}:%
\[
\lim_{\varepsilon\rightarrow0^{+}}f(\varepsilon)=\lim\limits_{s\rightarrow
1\cdot\infty}s\,F_{\text{alt}}(s).
\]
Differentiating, we obtain
\begin{align*}
f^{\prime}(x)  &  =\lambda(2-\rho)\left[  0-\theta(x)\right]  -\lambda
f(x)\left[  1-0\right]  -\lambda%
{\displaystyle\int\limits_{0}^{x}}
f(t)\left[  0-\theta(x-t)\right]  dt\\
&  =-\lambda(2-\rho)\theta(x)-\lambda f(x)+\lambda%
{\displaystyle\int\limits_{0}^{x}}
f(t)\theta(x-t)dt\\
&  =-\lambda(2-\rho)\theta(x)-\lambda f(x)+\lambda%
{\displaystyle\int\limits_{0}^{x}}
f(t)%
{\displaystyle\sum\limits_{k=1}^{\infty}}
e^{-k\,\rho}\dfrac{(k\,\rho)^{k-1}}{k!}\delta(x-t-a\,k)dt\\
&  =-\lambda(2-\rho)\theta(x)-\lambda f(x)+\lambda%
{\displaystyle\sum\limits_{k=1}^{\infty}}
e^{-k\,\rho}\dfrac{(k\,\rho)^{k-1}}{k!}%
{\displaystyle\int\limits_{0}^{x}}
f(t)\delta(x-t-a\,k)dt\\
&  =-\lambda(2-\rho)\theta(x)-\lambda f(x)+\lambda%
{\displaystyle\sum\limits_{k=1}^{\infty}}
e^{-k\,\rho}\dfrac{(k\,\rho)^{k-1}}{k!}f(x-a\,k).
\end{align*}
For $0<x<a$,%
\[%
\begin{array}
[c]{ccc}%
f^{\prime}(x)=-\lambda\,f(x), &  & f(0^{+})=\lambda
\end{array}
\]
implies%
\[
f(x)=\lambda\,e^{-\lambda\,x}.
\]
Note that $\lim_{\varepsilon\rightarrow0^{+}}f(a\,k+\varepsilon)=0$ for each
$k\geq1$ because, if a client arrives at the same moment the server becomes
available, the client is taken immediately (by LIFO) and there is no waiting.
For $a<x<2a$,%
\begin{align*}
f^{\prime}(x)  &  =-\lambda\,f(x)+\lambda\,e^{-\rho}\cdot\lambda
\,e^{-\lambda(x-a)}\\
&  =-\lambda\,f(x)+\lambda^{2}e^{-\lambda\,x}%
\end{align*}
coupled with $f(a^{+})=0$ implies%
\[
f(x)=\lambda^{2}(x-a)e^{-\lambda\,x}.
\]
For $2a<x<3a$,%
\begin{align*}
f^{\prime}(x)  &  =-\lambda\,f(x)+\lambda\,e^{-\rho}\cdot\lambda
^{2}(x-2a)e^{-\lambda(x-a)}+\lambda\,e^{-2\rho}\frac{2\rho}{2!}\cdot
\lambda\,e^{-\lambda(x-2a)}\\
&  =-\lambda\,f(x)+\lambda^{3}(x-a)e^{-\lambda\,x}%
\end{align*}
coupled with $f(2a^{+})=0$ implies%
\[
f(x)=\frac{1}{2}\lambda^{3}x(x-2a)e^{-\lambda\,x}.
\]
For $3a<x<4a$,%
\begin{align*}
f^{\prime}(x)  &  =-\lambda\,f(x)+\lambda\,e^{-\rho}\cdot\frac{1}{2}%
\lambda^{3}(x-a)(x-3a)e^{-\lambda(x-a)}\\
&  +\lambda\,e^{-2\rho}\frac{2\rho}{2!}\cdot\lambda^{2}(x-3a)e^{-\lambda
(x-2a)}+\lambda\,e^{-3\rho}\frac{(3\rho)^{2}}{3!}\cdot\lambda\,e^{-\lambda
(x-3a)}\\
&  =-\lambda\,f(x)+\frac{1}{2}\lambda^{4}x(x-2a)e^{-\lambda\,x}%
\end{align*}
coupled with $f(3a^{+})=0$ implies%
\[
f(x)=\frac{1}{6}\lambda^{4}x^{2}(x-3a)e^{-\lambda\,x}.
\]
More generally, for $k\,a<x<(k+1)a$, we obtain%
\[
f(x)=\frac{1}{k!}\lambda^{k+1}x^{k-1}(x-k\,a)e^{-\lambda\,x}%
\]
and thus the waiting time density for LIFO\ is completely understood.
\ Prabhu\ \cite{Pra1-que2} and Peters \cite{Pet-que2} evidently hold priority
in discovering this formula, the latter correcting an error in
\cite{Rio1-que2}. The density for M/D/1-FIFO is likewise completely understood
\cite{Erl-que2} -- no surprises occur here -- although it is M/U/1-FIFO
densities which are only \textit{effectively} known \cite{Fi-que2}.
\ Stitching the fragments together gives the LIFO density function pictured in
Figure 1, for parameter values $\lambda=2$ and $\mu=3$; hence $\rho=2/3$ and
$a=1/3$. \ It is interesting to compare this plot with Figure 3 of
\cite{Fi-que2}, the FIFO\ density for identical parameter values.

Let $\xi=a^{2}$ and $\eta=a^{3}$. \ Moments of $W_{\text{que}}$ for LIFO\ are
\cite{TS-que2, Coh-que2}\
\[%
\begin{array}
[c]{ccc}%
\text{mean}=-F^{\prime}(0)=\dfrac{\lambda\,\xi}{2(1-\rho)}, &  &
\text{variance}=F^{\prime\prime}(0)-F^{\prime}(0)^{2}=\dfrac{\lambda\,\eta
}{3(1-\rho)^{2}}+\dfrac{\lambda^{2}\xi^{2}(1+\rho)}{4(1-\rho)^{3}}%
\end{array}
\]
giving $\frac{1}{3}$ and $\frac{7}{9}$ respectively. \ The mean of
$W_{\text{que}}$ for FIFO is the same as that for LIFO; the corresponding
variance is smaller:%
\[%
\begin{array}
[c]{c}%
\dfrac{\lambda\,\eta}{3(1-\rho)}+\dfrac{\lambda^{2}\xi^{2}}{4(1-\rho)^{2}}%
\end{array}
\]
giving $\frac{5}{27}$. \ ($\xi$ and $\eta$ are second and third service time
moments, used for consistency with earlier work.)

We clarify that, while $\theta(x)$ is the notation for service time density in
\cite{Fi-que2}, $\theta(x)$ here is the notation for busy period length
density. \ The probability that \{a busy period $X$ is of length exactly
$k\,a$\} is equal to the probability that \{exactly $k$ clients of a queue,
having precisely $1$ starting client and traffic intensity $\rho$, are served
before the queue first vanishes\}. \ This is known as the Borel distribution
\cite{CG-que2, Pra2-que2, Tan-que2}, which satisfies \cite{Tak-que2, Pac-que2,
SH-que2}%
\[
\mathbb{E}\left(  \exp\left(  -s\,\frac{X}{a}\right)  \right)  =-\dfrac
{1}{\rho}\,\omega\left(  -\rho\exp\left(  -\rho-s\right)  \right)
\]
for M/D/1; thus the formula for $\Theta(s)=\mathbb{E}\left(  \exp\left(
-s\,X\right)  \right)  $ follows.

Study of busy period lengths is possible for M/U/1 -- the density involves
modified Bessel functions of the first kind \cite{ACW-que2} -- it is far more
complicated than the density for M/D/1.

\section{SIRO}

The probability density function $f(x)$ of $W_{\text{que}}$ has Laplace
transform \cite{Ki2-que2, LGa1-que2, LGa2-que2, BFL-que2}%
\[
F(s)=1+\frac{\lambda(1-\rho)}{s}%
{\displaystyle\int\limits_{\Theta(s)}^{1}}
\Phi(s,z)\Psi(s,z)dz=F_{\text{alt}}(s)+1-\rho
\]
where%
\[
\Phi(s,z)=\frac{1-z-a\,s/(1-\rho)}{z-e^{-a\,s-\rho(1-z)}}-\frac{1-z}%
{z-e^{-\rho(1-z)}},
\]%
\[
\Psi(s,z)=\exp\left(  -%
{\displaystyle\int\limits_{z}^{1}}
\frac{dy}{y-e^{-a\,s-\rho(1-y)}}\right)  .
\]
The integral underlying $\Psi(s,z)$ is intractable; our symbolic approach for
FIFO\ \&\ LIFO seems inapplicable for SIRO.

We therefore turn to a numeric approach developed by Burke \cite{Bur-que2},
which is based on certain simplifying properties of M/D/1 that do not easily
generalize. \ Since $\mu=1$ is assumed in \cite{Bur-que2}, we take
$\lambda=2/3$ and rescale at the end. \ Two recursions are key:%
\[
P_{n}=P_{n-1}e^{\lambda}-(P_{0}+P_{1})\frac{\lambda^{n-1}}{(n-1)!}-%
{\displaystyle\sum\limits_{j=2}^{n-1}}
P_{j}\frac{\lambda^{n-j}}{(n-j)!},
\]%
\[%
\begin{array}
[c]{ccc}%
P_{0}=1-\lambda, &  & P_{1}=P_{0}e^{\lambda}-P_{0}%
\end{array}
\]
where the empty sum convention holds for $n=2$, and
\[
Q_{i}(n,m)=\left\{
\begin{array}
[c]{lll}%
\left(  1-\dfrac{1}{n}\right)  e^{-\lambda}\,%
{\displaystyle\sum\limits_{k=0}^{m}}
\,\dfrac{\lambda^{k}}{k!}\,Q_{i-1}(n+k-1,m) &  & \text{if }n>0,\\
0 &  & \text{otherwise;}%
\end{array}
\right.
\]%
\[%
\begin{array}
[c]{ccc}%
Q_{i}(1,m)=\delta_{i}, &  & Q_{0}(n,m)=\left\{
\begin{array}
[c]{lll}%
\dfrac{1}{n} & \smallskip & \text{if }n>0,\\
0 &  & \text{otherwise}%
\end{array}
\right.
\end{array}
\]
where $\delta_{i}$ is the Kronecker delta and $m$ is suitably large. \ These
lead to the probability that waiting time is $\leq t$:%
\[
H(t,m)=(1-\lambda)+\lambda\,%
{\displaystyle\sum\limits_{n=1}^{\infty}}
P_{n-1}\cdot\left[  \left(  t-\left\lfloor t\right\rfloor \right)
Q_{\left\lfloor t\right\rfloor }(n,m)+%
{\displaystyle\sum\limits_{i=0}^{\left\lfloor t\right\rfloor -1}}
Q_{i}(n,m)\right]
\]
in the limit as $m\rightarrow\infty$, where again the empty sum convention
holds for $0\leq t<1$. \ 

Clearly the preceding expression simplifies to%
\[
(1-\lambda)+\lambda\,t\,%
{\displaystyle\sum\limits_{n=1}^{\infty}}
\frac{P_{n-1}}{n}%
\]
for all $0\leq t\leq1$ and any $m\geq0$. \ The dependence of $H(t,m)$ on $m$
becomes visible for $t>1$. \ To find the waiting time density value, let
$t^{\prime}=t+\frac{1}{4}$ (an arbitrary choice) and assume $0<t<t^{\prime}%
<1$. \ Then%
\[
4\left(  H(t^{\prime},m)-H(t,m)\right)  =\lambda\,%
{\displaystyle\sum\limits_{n=1}^{\infty}}
\frac{P_{n-1}}{n}%
\]
and, rescaling,%
\[
12\left(  H(t^{\prime},m)-H(t,m)\right)  =2\,%
{\displaystyle\sum\limits_{n=1}^{\infty}}
\frac{P_{n-1}}{n}=1.176773....
\]
This corresponds to the height of the leftmost rectangle in Figure 2,
surmounting the interval $[0,\frac{1}{3}]$.

For $1<t<t^{\prime}<2$, the rescaled density value is
\begin{align*}
\lim_{m\rightarrow\infty}\,12\left(  H(t^{\prime},m)-H(t,m)\right)   &  =2\,%
{\displaystyle\sum\limits_{n=2}^{\infty}}
P_{n-1}\left(  1-\dfrac{1}{n}\right)  e^{-\lambda}\,%
{\displaystyle\sum\limits_{k=0}^{\infty}}
\,\dfrac{\lambda^{k}}{k!}\,\frac{1}{n+k-1}\\
&  =0.392099...
\end{align*}
which corresponds to the height of the rectangle surmounting the interval
$[\frac{1}{3},\frac{2}{3}]$. \ Unlike earlier, the requirement that $m$ is
large becomes essential. \ For $2<t<t^{\prime}<3$, the formula for
$Q_{2}(n,m)$ involving $\{Q_{0}(n,m),Q_{1}(n,m)\}$ is more complicated than
that for $Q_{1}(n,m)$. \ It is better to implement the recursion than to write
the formula. \ Summarizing, initial segments of the sought-after density
function are
\[
g(x)=\left\{
\begin{array}
[c]{lll}%
c_{1}=1.176773... &  & \text{if }0<x<\frac{1}{3},\\
c_{2}=0.392099... &  & \text{if }\frac{1}{3}<x<\frac{2}{3},\\
c_{3}=0.179926... &  & \text{if }\frac{2}{3}<x<1,\\
c_{4}=0.095228... &  & \text{if }1<x<\frac{4}{3},\\
c_{5}=0.054829... &  & \text{if }\frac{4}{3}<x<\frac{5}{3},\\
c_{6}=0.033415... &  & \text{if }\frac{5}{3}<x<2,\\
c_{7}=0.021228... &  & \text{if }2<x<\frac{7}{3},\\
c_{8}=0.013923... &  & \text{if }\frac{7}{3}<x<\frac{8}{3},\\
c_{9}=0.009369... &  & \text{if }\frac{8}{3}<x<3.
\end{array}
\right.
\]
The Laplace transform of $g(x)$ is%
\[
G(s)=%
{\displaystyle\sum\limits_{j=1}^{\infty}}
c_{j}\frac{\exp\left(  -\frac{j-1}{3}s\right)  -\exp\left(  -\frac{j}%
{3}s\right)  }{s}%
\]
and verification that $G(s)=F_{\text{alt}}(s)$ can be done experimentally
(although not yet theoretically). \ 

We wonder about possible generalizations of this approach:\ works that cite
\cite{Bur-que2} include \cite{Lee-que2, CC-que2, DL-que2}. \ It is known (by
other techniques) that the mean of $W_{\text{que}}$ for SIRO is the same as
that for FIFO and LIFO; the corresponding variance is between the two extremes
\cite{Cn2-que2, Fuh-que2, TK-que2}:
\[
\dfrac{2\,\lambda\,\eta}{3(1-\rho)(2-\rho)}+\dfrac{\lambda^{2}\xi^{2}(2+\rho
)}{4(1-\rho)^{2}(2-\rho)}%
\]
giving $\frac{1}{3}$.%
\begin{figure}
[ptb]
\begin{center}
\includegraphics[
height=4.3145in,
width=6.634in
]%
{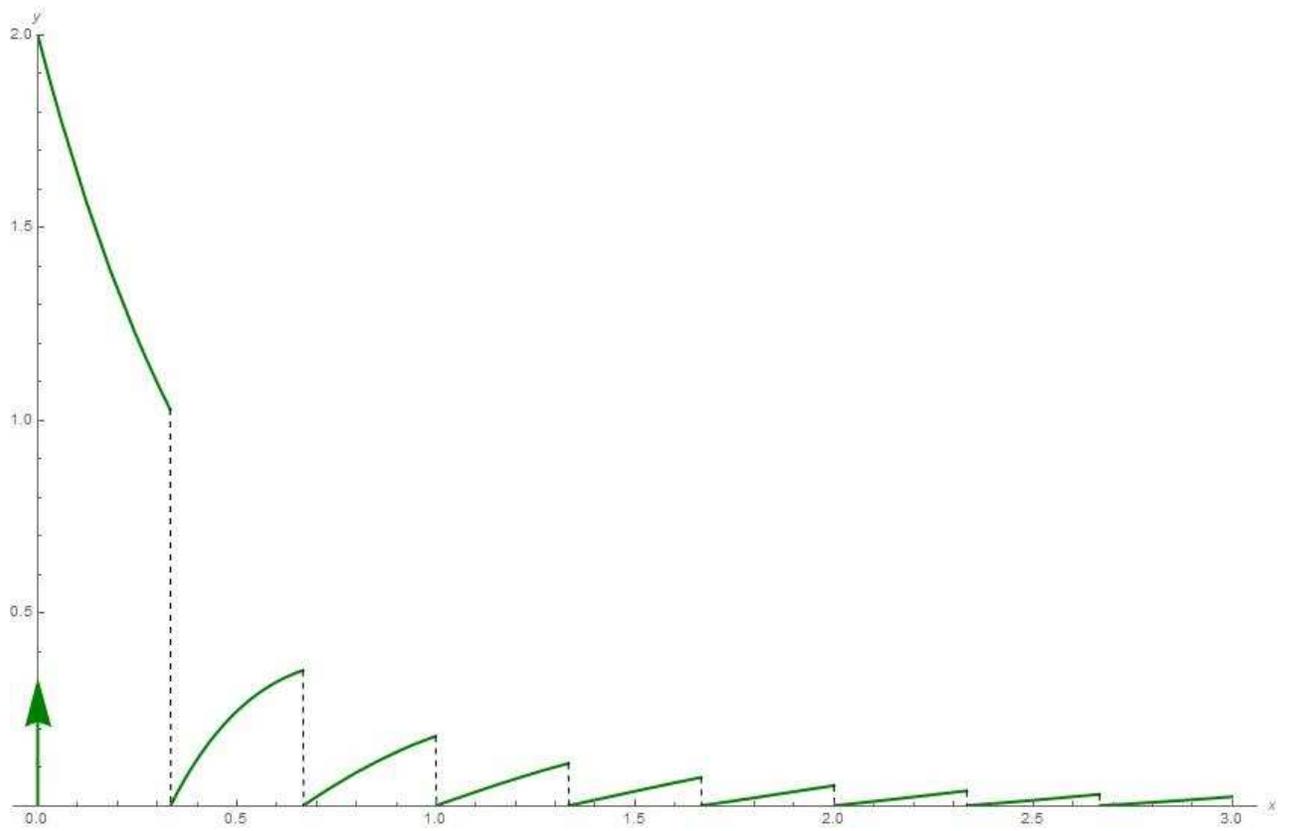}%
\caption{Waiting time density plot $y=f(x)$ for Deterministic[$\frac{1}{3}$]
last-in-first-out service.}%
\end{center}
\end{figure}
\begin{figure}
[ptb]
\begin{center}
\includegraphics[
height=4.3725in,
width=6.7239in
]%
{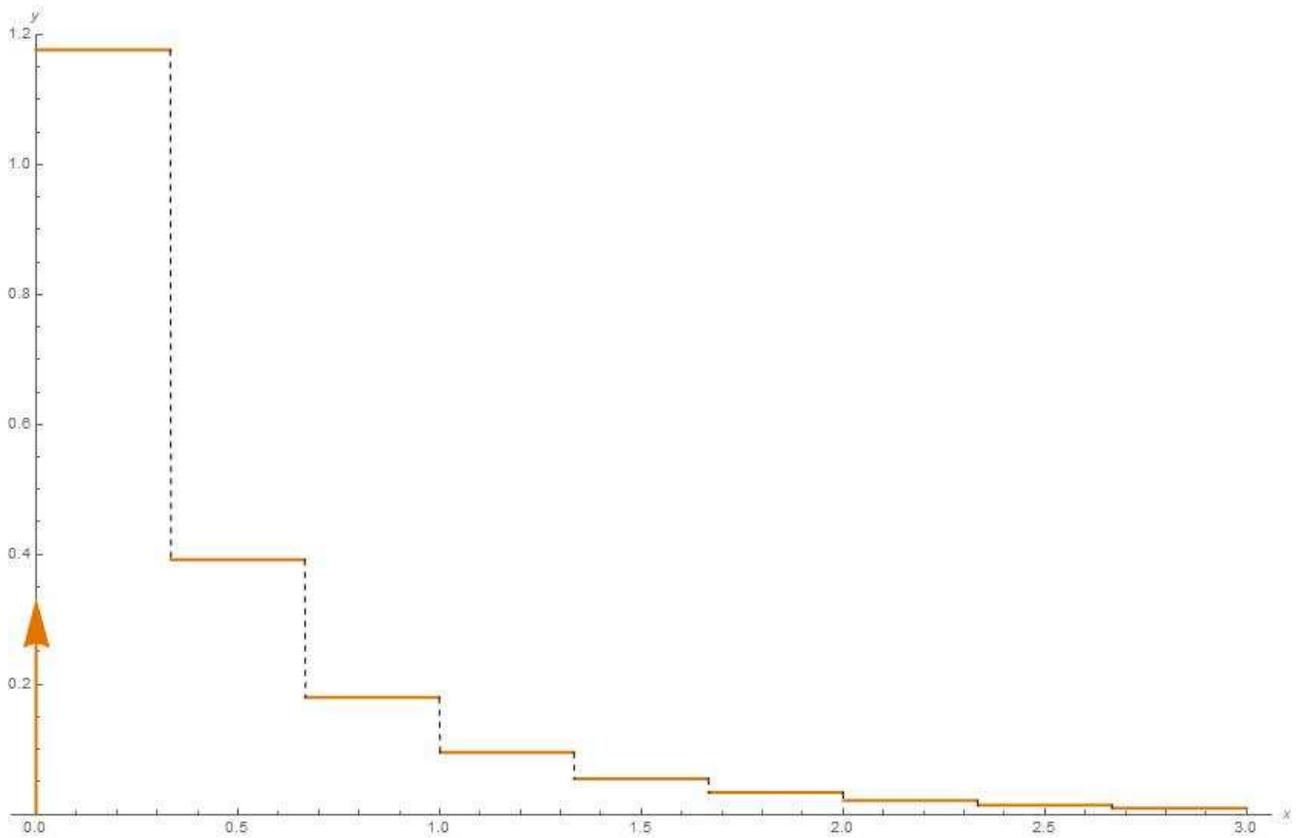}%
\caption{Waiting time density plot $y=g(x)$ for Deterministic[$\frac{1}{3}$]
serve-in-random-order policy.}%
\end{center}
\end{figure}

\section{Addendum}

Let us briefly summarize corresponding results for M/M/1. \ For LIFO, the
density $f(x)$ for $W_{\text{que}}$ is \cite{Rio1-que2, Vau-que2, Rio2-que2,
Kst-que2}%
\[
f(x)=(1-\rho)\delta(x)+\sqrt{\frac{\mu}{\lambda}}\,\frac{\exp\left[
-(\lambda+\mu)x\right]  \,I_{1}\left(  2\sqrt{\lambda\,\mu}\,x\right)  }{x}%
\]
where $I_{1}(\cdot)$ is the modified Bessel function of first order. \ Since
$\xi=2/\mu^{2}$ and $\eta=6/\mu^{3}$, the mean and variance are $\frac{2}{3}$
and $\frac{32}{9}$ respectively for parameter values $\lambda=2$ and $\mu=3$.
\ Note that $\frac{32}{9}>\frac{8}{9}$, which is the analogous variance for FIFO.

For SIRO, the density $f(x)$ for $W_{\text{que}}$ is likewise of the form
$(1-\rho)\delta(x)+\varphi(x)$ where \cite{Flat-que2, BMN-que2}%
\[
\varphi(x)=2(\mu-\lambda)%
{\displaystyle\int\limits_{0}^{\pi}}
\,\frac{\exp\left[  \left(  2\psi(\tau)-\tau\right)  \cot(\tau)\right]  }%
{\exp\left[  \pi\cot(\tau)\right]  +1}\,\frac{\exp\left[  -\lambda\left(
1-\frac{2}{\sqrt{\rho}}\cos(\tau)+\frac{1}{\rho}\right)  x\right]  }%
{1-\frac{2}{\sqrt{\rho}}\cos(\tau)+\frac{1}{\rho}}\sin(\tau)d\tau
\]
and%
\[
\psi(\tau)=\left\{
\begin{array}
[c]{lll}%
\arctan\left(  \dfrac{\sin(\tau)}{\cos(\tau)-\sqrt{\rho}}\right)  &  &
\text{if }0\leq\tau\leq\arccos\left(  \sqrt{\rho}\right)  ,\\
\pi+\arctan\left(  \dfrac{\sin(\tau)}{\cos(\tau)-\sqrt{\rho}}\right)  &  &
\text{if }\arccos\left(  \sqrt{\rho}\right)  <\tau\leq\pi.
\end{array}
\right.
\]
The mean is the same as that for LIFO; the variance is $\frac{14}{9}$ for
$\lambda=2$ and $\mu=3$, intermediate to $\frac{8}{9}$ and $\frac{32}{9}$.
\ Plotting the (monotone decreasing) densities for FIFO, LIFO and SIRO
together, it becomes evident that

\begin{itemize}
\item both short \&\ long waiting times occur more often under SIRO\ than
under FIFO

\item both short \&\ long waiting times occur more often under LIFO\ than
under SIRO
\end{itemize}

\noindent consistent with intuition.

We mention an unanswered question from \cite{Fn-que2}:\ for M/M/$c$-FIFO, what
is the maximum wait time density associated with the queue over a specified
time interval $[0,n]$? \ This open problem can be extended to LIFO\ and
SIRO\ as well. \ See also \cite{Fc-que2} for treatment (analogous to M/D/1
here)\ of the D/M/1 queue.

\section{Acknowledgements}

I am grateful to innumerable software developers. \ Mathematica routines
NDSolve for delay-differential equations and InverseLaplaceTransform (for Mma
version $\geq12.2$) plus ILTCME \cite{CME-que2} assisted in numerically
confirming many results. \ R steadfastly remains my favorite statistical
programming language.

\end{document}